\documentclass [12pt,a4paper] {article}
\usepackage[affil-it]{authblk}

\usepackage [T1] {fontenc}
\usepackage {amssymb}
\usepackage {amsmath}
\usepackage {amsthm}
\usepackage{hyperref}
\usepackage{tikz}

\usetikzlibrary{arrows,decorations.pathmorphing,backgrounds,positioning,fit,petri}

\usetikzlibrary{matrix}

\DeclareMathOperator {\td} {td}

\DeclareMathOperator {\ldim} {l.dim}

\DeclareMathOperator {\Exp} {Exp}

\DeclareMathOperator {\D} {D}

\DeclareMathOperator {\pr} {pr}

\DeclareMathOperator {\rank} {rank}
\DeclareMathOperator {\E} {E}

\theoremstyle {definition}
\newtheorem {definition}{Definition} [section]

\theoremstyle {plain}

\newtheorem {lemma} [definition] {Lemma}
\newtheorem {theorem} [definition] {Theorem}
\newtheorem {proposition} [definition] {Proposition}

\theoremstyle {remark}

\newtheorem* {notation} {Notation}

\newtheorem{innercustomthm}{Theorem}
\newenvironment{customthm}[1]
  {\renewcommand\theinnercustomthm{#1}\innercustomthm}
  {\endinnercustomthm}

\makeatletter
\newcommand {\forksym} {\raise0.2ex\hbox{\ooalign{\hidewidth$\vert$\hidewidth\cr\raise-0.9ex\hbox{$\smile$}}}}
\def\@forksym@#1#2{\mathrel{\mathop{\forksym}\displaylimits_{#2}}}
\def\forkind{\@ifnextchar_{\@forksym@}{\forksym}}

\makeatother

\makeatletter
\def\@maketitle{%
  \newpage
  \null
  \vskip 2em%
  \begin{center}%
  \let \footnote \thanks
    {\Large\bfseries \@title \par}%
    \vskip 1.5em%
    {\normalsize
      \lineskip .5em%
      \begin{tabular}[t]{c}%
        \@author
      \end{tabular}\par}%
    \vskip 1em%
    {\normalsize \@date}%
  \end{center}%
  \par
  \vskip 1.5em}
\makeatother

\begin {document}

\title{Ax-Schanuel for Linear Differential Equations}

\author{Vahagn Aslanyan\footnote{E-mail: \texttt{vahagn.aslanyan@gmail.com}}}
\affil{Mathematical Institute, University of Oxford, Oxford, OX2 6GG, UK\\
Institute of Mathematics of National Academy of Sciences of RA, Yerevan, Armenia
}

\maketitle

\begin{abstract}
We generalise the exponential Ax-Schanuel theorem to arbitrary linear differential equations with constant coefficients. Using the analysis of the exponential differential equation by J. Kirby (\cite{Kirby-thesis, Kirby-semiab}) and C. Crampin (\cite{Crampin}) we give a complete axiomatisation of the first order theories of linear differential equations and show that the generalised Ax-Schanuel inequalities are adequate for them.

\end{abstract}

\section{Introduction}\label{intro}
\setcounter{equation}{0}

In \cite{Lang-tr} Serge Lang mentions that Stephen Schanuel conjectured that for any $\mathbb{Q}$-linearly independent complex numbers $z_1,\ldots,z_n$ one has
\begin{equation}\label{Schanuels-conjecture}
\td_{\mathbb{Q}}\mathbb{Q}(z_1,\ldots,z_n,e^{z_1},\ldots,e^{z_n}) \geq n.
\end{equation}

This is now known as Schanuel's conjecture. It generalises many results and conjectures in transcendental number theory and is wide open.

Schanuel's conjecture (and its real version) is closely related to the model theory of the complex (real) exponential field $\mathbb{C}_{\exp}=(\mathbb{C};+,\cdot,\exp)$ (respectively, $\mathbb{R}_{\exp}=(\mathbb{R};+,\cdot,\exp)$, see \cite{Mac-Wil}). Most pertinent here is Boris Zilber's approach. He noticed that the inequality \eqref{Schanuels-conjecture} states the positivity of a predimension in the sense of Hrushovski. More precisely, Schanuel's conjecture is equivalent to the following statement: for any $z_1,\ldots,z_n \in \mathbb{C}$ the inequality
\begin{equation}\label{Schanuel-predimension}
\delta(\bar{z}) := \td_{\mathbb{Q}}\mathbb{Q}(\bar{z},\exp(\bar{z}))-\ldim_{\mathbb{Q}}(\bar{z}) \geq 0
\end{equation}
holds, where $\td$ and $\ldim$ stand for transcendence degree and linear dimension respectively. Here $\delta$ satisfies the submodularity law which allows one to carry out a Hrushovski construction. In this way B. Zilber constructed pseudo-exponentiation on algebraically closed fields of characteristic zero. He proved that there is a unique model of that (not first-order) theory in each uncountable cardinality and conjectured that the model of cardinality $2^{\aleph_0}$ is isomorphic to $\mathbb{C}_{\exp}$. Since \eqref{Schanuel-predimension} holds for pseudo-exponentiation (it is included in the axiomatisation given by Zilber), Zilber's conjecture implies Schanuel's conjecture. For details on pseudo-exponentiation see \cite{Zilb1,Zilb2,Zilb3,Zilb4}.

Though Schanuel's conjecture seems to be out of reach of modern methods in mathematics, James Ax proved its differential analogue in 1971 (\cite{Ax}). It is now known as the Ax-Schanuel theorem or inequality.

\begin{theorem}[Ax-Schanuel]
Let $\mathcal{K}=(K;+,\cdot,\D)$ be a differential field with field of constants $C$. If $(x_1,y_1),\ldots,(x_n,y_n)$ are non-constant solutions to the exponential differential equation $\D y = y \D x$ then
\begin{equation}\label{Ax-Schanuel}
\td_CC(\bar{x},\bar{y})-\ldim_{\mathbb{Q}}(\bar{x}/C) \geq 1,
\end{equation}
where $\ldim_{\mathbb{Q}}(\bar{x}/C)$ is the dimension of the $\mathbb{Q}$-span of $x_1, \ldots, x_n$ in the quotient vector space $K/C$.
\end{theorem}

Here again we have a predimension inequality, which will be part of the first order theory of the reduct $\mathcal{K}_{\Exp}=(K;+,\cdot,\Exp)$ of $\mathcal{K}$ where $\Exp$ is a binary predicate for the set of solutions of the exponential differential equation, i.e. it is interpreted in a differential field $\mathcal{K}$ as
$$\Exp(K)=\{ (x,y)\in K^2: \D y = y\D x \}.$$ 

Therefore a natural question arises: if one carries out a Hrushovski construction with this predimension and class of reducts, will one end up with a similar reduct of a differentially closed field?\footnote{{This question arises from similarity to the pseudo-exponentiation.}} Zilber calls predimensions with this property \textit{adequate}. 

\begin{definition}[\cite{Aslanyan-thesis}]\label{def-adequate}
{A predimension inequality on the set of solutions of a differential equation $E(x,y)$ is \textit{adequate}, if the models of the theory of $E$-reducts $(K;+,\cdot, E)$ of differential fields $(K;+,\cdot, \D)$ satisfy the strong amalgamation property, and the Hrushovski construction yields an $E$-reduct of a differentially closed field, that is, the Hrushovski-Fra\"{i}ss\'{e} limit $\mathcal{U}$ is isomorphic to such a reduct.\footnote{{This corresponds to \textit{strong adequacy} in \cite{Aslanyan-thesis}, and there adequacy means that $\mathcal{U}$ is elementarily equivalent to an $E$-reduct of a differentially closed field. However, in this paper we do not need this weaker notion and prefer the term ``adequate'' to ``strongly adequate'' for simplicity.}}}
\end{definition}


Thus the question is whether the Ax-Schanuel inequality is adequate. {Cecily Crampin studied the exponential differential equation in her PhD thesis \cite{Crampin} and gave a criterion for a system of exponential differential equations to have a solution (analogous to pseudo-exponentiation), known as \emph{existential} or \emph{exponential closedness} (in fact, it is a special case of the full existential closedness property proved by Kirby). She also considered the predimension function and proved some form of strong existential closedness for reducts of differentially closed fields. Her results can be used to show that Ax-Schanuel is adequate, though it is not explicit in her work as she does not construct the strong Fra\"iss\'e limit.}

Jonathan Kirby considered the same problem in a much more general context. He studied exponential differential equations for semiabelian varieties, observed that Ax-Schanuel holds in this setting too and proved that it is adequate along with giving an axiomatisation of the complete theory of the reducts (see \cite{Kirby-thesis, Kirby-semiab}). The axiomatisation is again very similar to pseudo-exponentiation. An important property that implies adequacy of Ax-Schanuel is the \textit{strong existential closedness} which means that models of $T_{\Exp}$ are existentially closed in strong extensions {(see Section \ref{Predimension} for the definition of strong extensions)}. More details on this, in particular an axiomatisation of $T_{\Exp}$, will be presented in Section \ref{exp}.

Once this is done, one naturally asks the question of whether something similar can be done for other differential equations.\footnote{{This is the main question of my PhD thesis posed by my supervisors Boris Zilber and Jonathan Pila. See \cite{Aslanyan-thesis} for more details.}}  In other words, one wants to find adequate predimension inequalities for differential equations. In this paper we show that this can be done for any linear differential equations with constant coefficients (the exponential differential equation being a special case of it). 

We formulate our main results below. For a differential field $\mathcal{K}$ and a non-constant element $x \in K$ define a derivation $\partial_x : K \rightarrow K$ by $\partial_x = (\D x)^{-1} \cdot \D$. Then consider the differential equation
\begin{equation}\label{eqn}
(\D x)^{2n-1} \left[ \partial_x^ny +c_{n-1}\partial_x^{n-1}y+\ldots+c_1\partial_xy+c_0y \right]=0,\tag{*}
\end{equation}
where the coefficients are constants with $c_0 \neq 0$.
The role of the factor $(\D x)^{2n-1}$ is that it clears the denominator of the left hand side of
\eqref{eqn}, making it a differential polynomial. Let $\lambda_1,\ldots, \lambda_n$ be the roots of the characteristic polynomial $p(\lambda)=\lambda^n+\sum_{0\leq i < n} c_i \lambda^i$. Then the Ax-Schanuel theorem for \eqref{eqn} is as follows.

\begin{customthm}{\ref{Ax-Schanuel-for-higher-order}}\label{Ax-intro}
Let $(x_i,y_i),~ i=1,\ldots,m,$ be non-constant solutions to the equation \eqref{eqn} in a differential field $\mathcal{K}$ such that $y_i, \partial_{x_i} y_i, \ldots, \partial_{x_i}^{n-1}y_i$ are linearly independent over $C$ for every $i$. Then
\begin{equation}\label{inequality-intro}
\td_CC(\bar{x},\bar{y},\partial_{\bar{x}}\bar{y},\ldots,\partial_{\bar{x}}^{n-1}\bar{y}) \geq \ldim_{\mathbb{Q}}(\lambda_1 \bar{x},\ldots,\lambda_n \bar{x}/C)+1, \tag{**}
\end{equation}
where $\partial^j_{\bar{x}}\bar{y}:=(\partial^j_{x_1}y_1, \ldots, \partial^j_{x_m}y_m)$.
\end{customthm}

Our problem is to axiomatise the complete theory of the equation \eqref{eqn}. As the inequality \eqref{inequality-intro} suggests, it will be natural to work with reducts of differential fields with an $(n+1)$-ary relation $\E_n(x,y, y_1,\ldots,y_{n-1})$ {for the set of all $(n+1)$-tuples $(x,y,y_1,\ldots,y_{n-1})$ where $(x,y)$ is a solution of \eqref{eqn} and $y_i = \partial_x^i y$}. Then \eqref{inequality-intro} can be written in a first-order way. For the complete axiomatisation we will need an important axiom scheme called \textit{Existential Closedness}. For a structure $F$ in our language it can be formulated as follows.\footnote{For the definition of $\E_n$-$\Exp$-rotundity see Section \ref{Existential-closedness}. {Let us just mention here that it is an analogue of \textit{rotundity} which is called $\Exp$-rotundity in this paper. $\Exp$-Rotundity (first defined by Zilber under the name \textit{normality}) of a variety means that the variety is not ``too small'' in terms of the algebraic dimensions of its projections and some ``twisted'' projections. It is a sufficient condition for a variety to have a non-empty intersection with $\Exp$ in a differentially closed field. The appropriate conditions are dictated by the functional equation of the exponential function, namely $\exp(x+y)=\exp(x)\cdot \exp(y)$, and by Ax-Schanuel. Note also that $\Exp$-rotundity and $\E_n$-$\Exp$-rotundity are first-order properties in the language of rings.}}

\begin{itemize}
\item[EC'] For each $\E_n$-$\Exp$-rotund variety $ V \subseteq F^{m(n+1)}$ the intersection $V(F) \cap \E^m_n(F)$ is non-empty.
\end{itemize}

We will see that this axiom scheme along with the inequality \eqref{inequality-intro} and some basic axioms (which reveal the relationship between $\E_n$ and $\Exp$) will axiomatise the first order theory $T_{\E_n}$ of reducts of differentially closed fields in the corresponding language.

Our results rely heavily on the aforementioned analysis of the equation $\D y = y \D x$. In particular, we use the Ax-Schanuel theorem to prove Theorem \ref{Ax-intro}. Then we use the axiomatisation of $T_{\Exp}$ to obtain an axiomatisation of $T_{\E_n}$. In fact the exponential differential equation can be defined in our reducts which means that one can simply translate the axioms of $T_{\Exp}$ to obtain an axiomatisation of $T_{\E_n}$. However we give an axiomatisation based on the predimension inequality \eqref{inequality-intro} and show the adequacy of that predimension, which means that the $\E_n$-reduct of a countable saturated differentially closed field can be constructed by a Hrushovski construction with that predimension.

Let us also note that the problem of finding adequate predimension inequalities proved to be closely related to the problem of definability of the derivation of our differential field in its reducts. We have studied this question in \cite{Aslanyan}. In particular, the results of the current paper, with the analysis carried out in \cite{Aslanyan}, will show that $\D$ is not definable from linear differential equations.

{Finally, we remark that Kirby's work on the exponential differential equation is the functional analogue of Zilber's work on Schanuel's conjecture and complex exponentiation. Furthermore, Ax-Schanuel can be applied to deduce a weak version of Schanuel's conjecture. Namely, Kirby proved in \cite{Kirby-Schanuel} that Schanuel's conjecture has at most countable many ``essential'' counterexamples. This being said, it would be natural to ask if we can connect our results on linear differential equations with some number theoretic problems. The answer is that complex solutions of linear differential equations with constant coefficients are just linear combinations of exponential functions. Thus, we do not get any ``new'' functions and, therefore, it is unlikely that our theorems will yield a transcendence result beyond the original Ax-Schanuel.}\\\\

\textbf{Acknowledgements.}
I am very grateful to my supervisors Boris Zilber and Jonathan Pila for many helpful discussions.\\
This research was supported by the University of Oxford Dulverton Scholarship.

\section{The exponential differential equation}\label{exp}
\setcounter{equation}{0}
In this section we give an axiomatisation of the theory of the exponential differential equation. We will work in the language $\mathfrak{L}_{\Exp}:=\{+,\cdot,0,1,\Exp\}$ where $\Exp$ is a binary predicate which will be interpreted in a differential field $\mathcal{K}=(K;+,\cdot,0,1,\D)$ as the set $\{ (x,y)\in K^2: \D y = y \D x \}$. In this case the reduct of $\mathcal{K}$ to the language $\mathfrak{L}_{\Exp}$ will be denoted by $\mathcal{K}_{\Exp}$. For a differentially closed field $\mathcal{K}$ we denote the complete first-order theory of $\mathcal{K}_{\Exp}$ by $T_{\Exp}$. For an $\mathfrak{L}_{\Exp}$-structure $\mathcal{F}_{\Exp}$\footnote{We will normally add a subscript $\Exp$ in the notations of $\mathfrak{L}_{\Exp}$-structures to emphasise the fact that they are $\mathfrak{L}_{\Exp}$-structures and to distinguish them from $\mathfrak{L}_{\E_n}$-structures considered later. It does not mean that they are reducts of some differential fields unless we explicitly state that they are.} and for a natural number $n$ we let
$$\Exp^n(F):=\{ (\bar{x},\bar{y})\in F^{2n}: \mathcal{F}_{\Exp} \models \Exp(x_i,y_i) \text{ for each } i \}.$$


As we already mentioned, an axiomatisation of $T_{\Exp}$ has been given by Kirby \cite{Kirby-thesis,Kirby-semiab} and partially by Crampin \cite{Crampin} (Kirby's work is much more general, he studies exponential differential equations of semiabelian varieties). The original idea of such an axiomatisation is due to Zilber in the context of pseudo-exponentiation \cite{Zilb1}. We refer the reader to \cite{Kirby-semiab, Bays-Kirby-exp, Zilb1, Crampin} for details and proofs of the results presented in this section.

Throughout the paper $\mathcal{K}=(K;+,\cdot,\D, 0,1)$ will be a differential field and $C$ will denote its field of constants.

\begin{theorem}[\cite{Ax}, Theorem 3]\label{Ax-version1}
For any $x_i, y_i \in K,~ i=1,\ldots,n$, if $\mathcal{K} \models \bigwedge_{i=1}^n\Exp(x_i,y_i)$ and $\td_CC(\bar{x},\bar{y}) \leq n$ then there are integers $m_1,\ldots,m_n$, not all of them zero, such that $m_1x_1+\ldots+m_nx_n \in C$ or, equivalently, $y_1^{m_1}\cdot \ldots \cdot y_n^{m_n} \in C$.
\end{theorem}

This can be given a geometric formulation. For a field we let $\mathbb{G}_a$ be its additive group and $\mathbb{G}_m$ be the multiplicative group. Also for a natural number $n$ we denote $G_n:=\mathbb{G}_a^n \times \mathbb{G}_m^n$. Observe that for a differential field $\mathcal{K}$ the set $\Exp(K) \subseteq K^2$ is a subgroup of $G_2(K)$. Notice that $\prod y_i^{m_i} = c \in C$ means that $(y_1,\ldots,y_n)$ lies in a $C$-coset of the subgroup of $\mathbb{G}_m^n(K)$ defined by $\prod y_i^{m_i} = 1$. {Similarly, $\sum m_ix_i=c \in C$ means that $(x_1,\ldots,x_n)$ lies in a $C$-coset of the subgroup of $\mathbb{G}_a^n(K)$ defined by $\sum m_ix_i=0$.}

The tangent space of $\mathbb{G}_m^n$ at the identity can be identified with $\mathbb{G}_a^n$. For an algebraic subgroup $H$ of $\mathbb{G}_m^n$ its tangent space at the identity, denoted $T_eH$, is an algebraic subgroup of $\mathbb{G}_a^n$. 
The tangent bundle of $H$ will be denoted by $TH$. Also, for an integer $n$ we let $\Exp^n(K):=\{ (\bar{x},\bar{y})\in K^{2n}: \mathcal{K} \models \bigwedge_{i=1}^n\Exp(x_i,y_i) \}$.

These observations allow one to reformulate Theorem \ref{Ax-version1} in a geometric language.

\begin{theorem}[Ax-Schanuel - version 2]\label{Ax-version2}
Let $V \subseteq G_n(K)$ be an algebraic variety defined over $C$ with $\dim(V) \leq n$. If $(\bar{x},\bar{y}) \in V(K)\cap \Exp^n(K)$ then there is a proper algebraic subgroup $H$ of $\mathbb{G}_m^n$ such that $(\bar{x},\bar{y})$ lies in a $C$-coset of $TH$, that is, $\bar{y} \in \gamma H$ and $\bar{x} \in \gamma'+T_e H$ for some constant points $\gamma \in \mathbb{G}_m^n(C)$ and $\gamma' \in \mathbb{G}_a^n(C)$.
\end{theorem}

If $V$ is a variety as above and $V(K)\cap \Exp^n(K) \neq \emptyset$ then we say $V$ has an exponential point. The Ax-Schanuel theorem can be thought of as a necessary condition for a variety to have an exponential point. We will shortly present the existential closedness statement, which is a sufficient condition for this. But for now we consider some basic axioms for an $\mathfrak{L}_{\Exp}$-structure $\mathcal{F}_{\Exp}$.

\begin{itemize}
\item[A1] $F$ is an algebraically closed field of characteristic $0$.

\item[A2] $C:=C_F= \{ c \in F: \mathcal{F}_{\Exp} \models \Exp(c,1) \}$ is an algebraically closed subfield of $F$.

\item[A3] $\Exp(F) = \{ (x,y)\in F^2: \Exp(x,y)\}$ is a subgroup of $G_1(F)$ containing $G_1(C)$.

\item[A4] The fibres of $\Exp$ in $\mathbb{G}_a(F)$ and $\mathbb{G}_m(F)$ are cosets of the subgroups $\mathbb{G}_a(C)$ and $\mathbb{G}_m(C)$ respectively.

\item[AS] For any $x_i, y_i \in F,~ i=1,\ldots,n$, if $\mathcal{F}_{\Exp} \models \bigwedge_{i=1}^n\Exp(x_i,y_i)$ and $\td_CC(\bar{x},\bar{y}) \leq n$ then there are integers $m_1,\ldots,m_n$, not all of them zero, such that $m_1x_1+\ldots+m_nx_n \in C$.

\item[{NT}] {$F \supsetneq C$.}
\end{itemize}

{Note that AS can be given by an axiom scheme. A compactness argument gives a uniform version of the Ax-Schanuel theorem in differential fields. That is, given a parametric family of varieties $V(\bar{c})$ over $C$, there is a finite number $N$, such that if for some $\bar{c}$ we have $(\bar{x},\bar{y})\in V(\bar{c})$ and $\dim V(\bar{c})\leq n$ then $m_1x_1+\ldots+m_nx_n \in C$ for some integers $m_i$ with $|m_i| \leq N$. Thus, AS above should be understood as this uniform version of Ax-Schanuel.}

{These axioms basically constitute the universal part of $T_{\Exp}$ with the exception that A1 is $\forall \exists$ and NT is existential. Models of the theory A1-A4,AS will be called $\Exp$-\emph{fields}.}


Now we turn to existential closedness. For a $k \times n$ matrix $M$ of integers we define $[M]:G_n(F) \rightarrow G_k(F)$ to be the map given by $[M]:(\bar{x},\bar{y}) \mapsto (u_1,\ldots,u_k, v_1,\ldots, v_k)$ where
$$u_i = \sum_{j=1}^n m_{ij}x_j \mbox{ and } v_i = \prod_{j=1}^n y_j^{m_{ij}}.$$

\begin{definition}\label{rotund}
An irreducible variety $V \subseteq G_n(F)$ is {$\Exp$-\emph{rotund}} if for any $1 \leq k \leq n$ and any $k\times n$ matrix $M$ of integers $\dim [M](V) \geq \rank M$. If for any non-zero $M$ the stronger inequality $\dim [M](V) \geq \rank M + 1$ holds then we say $V$ is \emph{strongly $\Exp$-rotund}.
\end{definition}

The definition of $\Exp$-rotundity is originally due to Zilber though he initially used the word \textit{normal} for these varieties \cite{Zilb1}. The term \emph{rotund} was coined by Kirby in \cite{Kirby-semiab}. 

Strong $\Exp$-rotundity fits with the Ax-Schanuel inequality in the sense that it is a sufficient condition for a variety defined over $C$ to contain a non-constant exponential point. More precisely, if $F$ is differentially closed and $V \subseteq G_n(F)$ is a strongly $\Exp$-rotund variety defined over the constants, then the intersection $V(F) \cap \Exp^n(F)$ contains a non-constant point. {This is a form of an existential closedness statement, and the existential closedness property that we will use for the axiomatisation of $T_{\Exp}$ generalises this considering varieties that are not necessarily defined over $C$.}

The existential closedness property for an ${\Exp}$-field $\mathcal{F}_{\Exp}$ is as follows.

\begin{itemize}
\item[EC] For each irreducible $\Exp$-rotund variety $V \subseteq G_n(F)$ the intersection $V(F) \cap \Exp^n(F)$ is non-empty.
\end{itemize}
As noted above, $V$ is not necessarily defined over $C$ and the point in the intersection may be constant. 


$\Exp$-rotundity of a variety is a definable property {in the language of rings}. This allows one to axiomatise the above statement by a first-order axiom scheme. Reducts of differentially closed fields satisfy EC and it gives a complete theory {(in the language of reduct $\mathfrak{L}_{\Exp}=\{+,\cdot,0, 1, \Exp \}$)} together with the axioms mentioned above.

\begin{theorem}[\cite{Kirby-semiab}]\label{exp-axioms}
The theory $T_{\Exp}$ is axiomatised by the following axioms and axiom schemes: \emph{A1-A4, AS, EC, NT}.
\end{theorem}

We also define free varieties (\cite{Kirby-semiab, Zilb1}) and present a result from \cite{Kirby-semiab} below. Although this is not essential for our main results, we will use it to establish a similar fact for linear differential equations of higher order (Section \ref{Rotundity and freeness}) {which nevertheless contributes to our understanding of the general picture. Indeed, for free varieties $\Exp$-rotundity is a necessary condition for having an $\Exp$-point, and we will see that a similar fact holds for higher order linear equations too.}

\begin{definition}
An irreducible variety $V \subseteq G_n(K)$ (defined over $C$) is $\Exp$-\emph{free} if it does not have a generic (over $C$) point $(\bar{a}, \bar{b})$ for which
$$\sum m_i a_i \in C \mbox{ or } \prod y_i^{k_i} \in C$$
for some integers $m_i$ and $k_i$ (not all of them zero).
\end{definition}
Note that this notion corresponds to \emph{absolute freeness} in \cite{Kirby-thesis,Kirby-semiab}.

\begin{proposition}[\cite{Kirby-thesis,Kirby-semiab}]\label{exp-free}
Let $V$ be an $\Exp$-free variety defined over $C$. If $V$ has a generic (over $C$) exponential point then it is strongly $\Exp$-rotund.
\end{proposition}

Finally let us make an easy observation which will be useful later.

\begin{lemma}\label{exp-lemma}
Let $\mathcal{K}$ be a differentially closed field. If $V \subseteq G_n(K)$ is $\Exp$-rotund then for any constant $c \in C^{\times}$ there is a point $(\bar{a}, \bar{b}) \in V(K)$ such that $\mathcal{K}_{\Exp} \models \Exp(ca_i, b_i)$ for all $i$.
\end{lemma}
\begin{proof}
Let $L: K^{2n} \rightarrow K^{2n}$ be the map $(\bar{x}, \bar{y}) \mapsto (c\bar{x}, \bar{y})$. It is easy to check that $V':=L(V)$ is Zariski closed and $\Exp$-rotund. Therefore there is a point $(\bar{u}, \bar{v}) \in V'(K) \cap \Exp^n(K)$. If $a_i = c^{-1}u_i,~ b_i=v_i$, then $(\bar{a}, \bar{b}) \in V(K)$ and $\Exp(ca_i, b_i)$ holds.
\end{proof}

\section{Higher order linear differential equations}\label{Higher-order-section}
\setcounter{equation}{0}

In this section we will use some facts and notions from the theory of abstract linear differential equations in differential fields (see \cite{Mar-dif}, Section 4).

Let us start with a motivating example which will make it clear which differential equations we will consider. If $x(t)$ and $y(t)$ are complex analytic functions with $y(t)=\exp(x(t))$ then they satisfy the differential equation $\frac{d}{d t} y(t) = y(t) \cdot \frac{d}{d t} x(t)$. Since we are interested in non-constant solutions, this equation can be written as $\frac{dy}{dx} = \frac{\frac{d}{d t} y}{\frac{d}{d t} x} = y$. Now if we replace $\frac{d}{d t}$ with $\D$,\footnote{{Recall that we work in a differential field $\mathcal{K}=(K;+,\cdot, \D)$ with a single derivation $\D$ and field of constants $C$.}} we will obtain the abstract exponential differential equation $\frac{\D y}{\D x} = y$. Here we could also argue as follows. In the differential equation $\frac{dy}{dx} = y$ replace differentiation with respect to $x$, that is, $\frac{d}{dx}$ with $\frac{1}{\D x} \cdot \D$ to get $\frac{\D y}{\D x} = y$. If $x \in K$ is a non-constant element then $\partial_x = \frac{1}{\D x} \cdot \D$ is a derivation of $K$ and the exponential differential equation can be written as $\partial_x(y)=y$. Here $\partial_x$ can be thought of as abstract differentiation with respect to $x$. {Observe also that for each $x \in K\setminus C$ the field of constants with respect to $\partial_x$ coincides with the field of constants with respect to the original derivation $\D$, which is denoted by $C$.}


Now we want to generalise this to higher order linear differential equations with constant coefficients. Consider the equation
\begin{equation}
\frac{d^ny}{dx^n} +c_{n-1}\frac{d^{n-1}y}{dx^{n-1}}+\ldots+c_1\frac{dy}{dx}+c_0y=0.
\end{equation}

Its solutions are {complex} linear combinations of exponential functions. We want to form the corresponding abstract differential equations whose solutions will be analogues of those combinations. As above we replace $\frac{d}{dx}$ by $\partial_x$ to obtain the equation
$$\partial_x^ny +c_{n-1}\partial_x^{n-1}y+\ldots+c_1\partial_xy+c_0y=0.$$ 
The left hand side of this equation is a differential rational function with denominator $(\D x)^{2n-1}$. We multiply through by this factor to make the left hand side into a polynomial. {This will also be useful later when we define the field of constants in the appropriate reduct.} Thus we consider the abstract differential equation
\begin{equation}\label{main-eq}
\Delta(x,y):=(\D x)^{2n-1}\left[ \partial_x^ny +c_{n-1}\partial_x^{n-1}y+\ldots+c_1\partial_xy+c_0y \right]=0
\end{equation}
in a differential field $\mathcal{K}$. {We assume the coefficients $c_i \in C$ are constants.}

The notation $\partial_x$ may misleadingly suggest that $x$ is fixed in the equation \eqref{main-eq} which is not the case. It should be considered as a two-variable equation. We prefer this way of writing our equation since otherwise it would be cumbersome. Note however that $\Delta(x,y)$ is not linear as a two-variable differential polynomial, it is linear with respect to $y$ only. We will assume that $c_0 \neq 0$ in order to avoid any possible degeneracies (like $\D y =0$).


Observe that by introducing new variables $z_0,\ldots,z_n$ we can write \eqref{main-eq} as the following system of equations
\begin{equation}\label{main-eq2}
\begin{cases}
z_n+c_{n-1}z_{n-1}+ \ldots +c_1z_1+c_0z_0=0, \\
z_0=y,\\
\D z_{i} = z_{i+1} \D x,~ i=0,\ldots ,n-1.
\end{cases}
\end{equation}


Let $p(\lambda) = \lambda^n+c_{n-1}\lambda^{n-1}+\ldots+c_0$ be the characteristic polynomial of \eqref{main-eq}. Let $\lambda_1,\ldots, \lambda_n$ be its roots and let $\mu_1, \ldots, \mu_k$ be its distinct roots with multiplicities $n_1, \ldots, n_k$ respectively. Since we have assumed $c_0$ is non-zero, $\lambda_i$'s are also non-zero.


Now we establish some auxiliary results which will be used in the proof of the Ax-Schanuel theorem for the equation \eqref{main-eq}. Since it is a universal statement, we can assume without loss of generality that $\mathcal{K}$ is differentially closed. This is not very important but makes our arguments easier as we do not have to worry about the existence of solutions of differential equations.

\begin{lemma}\label{solutions}
Let $x$ be a non-constant element of $K$ and let $y_i \in K\setminus \{0\}$ be such that $\partial_x y_i = \mu_i y_i $ for $i=1,\ldots, k$. Then $\bigcup_{i=1}^k \{y_i,xy_i, \ldots, x^{n_i-1}y_i \}$ forms a fundamental system of solutions\footnote{This means that those solutions form a $C$-linear basis for the space of all solutions.} to $\Delta(x,y)=0$.
\end{lemma}

Though the proof is very similar to that in the complex setting (see, for example, \cite{Birkhoff}), we nevertheless present it here for completeness.
\begin{proof}

Since $x$ is non-constant, the equation \eqref{main-eq} can be written as $p(\partial_x)y=0$. The operator $p(\partial_x)$ can be factored as
$$p(\partial_x) = \prod_{i=1}^k (\partial_x - \mu_i)^{n_i}.$$
It is easy to see that for any $0 \leq l < n_i$
$$(\partial_x-\mu_i)^{n_i} (x^ly_i) = 0.$$
Hence we have $p(\partial_x)(x^l y_i)=0$ and thus we have $n$ solutions to $\Delta(x,y) = 0$. Now we prove they are linearly independent.

Assume 
\begin{equation}\label{independent-solutions}
\sum_{i=1}^k\sum_{j=0}^{n_i-1} a_{ij}x^jy_i =0
\end{equation}
for some constants $a_{ij}$. Let $i$ be such that there is a non-zero coefficient $a_{ij}$. Let $t$ be the biggest number with $a_{it}\neq 0$. Consider the operator 
$$q(\partial_x) = (\partial_x - \mu_i)^t \prod_{s \neq i} (\partial_x - \mu_s)^{n_s}.$$

Clearly 
$$q(\partial_x)(x^jy_r) = 
\begin{cases}
0, \mbox{ if } r \neq i \mbox{ or } j<t,\\
t! \cdot \prod_{s \neq i} (\mu_i - \mu_s)^{n_s}\cdot y_i \neq 0 , \mbox{ if } r=i,~ j=t.
\end{cases}
$$

Now applying $q(\partial_x)$ to \eqref{independent-solutions} we get $a_{it}=0$, a contradiction.

{Thus, we found $n$ linearly independent solutions. Since the order of the equation is $n$, the set of solutions is an $n$-dimensional $C$-vector space, therefore the above solutions form a basis for that vector space.}
\end{proof}

If $y_1, \ldots, y_k$ are as in Lemma \ref{solutions}, then for any non-zero constants $a_{ij}$ the set $\bigcup_{i=1}^k \{a_{i0}y_i,a_{i1}xy_i, \ldots, a_{i,n_i-1}x^{n_i-1}y_i \}$ is a fundamental system of solutions to our equation. This kind of fundamental systems will be called \emph{canonical}. There is a unique such system up to multiplication by constants. Note also that we will treat canonical fundamental systems as ordered tuples, rather than as sets. Thus if we say $v_1, \ldots, v_n$ is a canonical fundamental system, then we mean that the first $n_1$ elements coincide (up to constants) with $y_1,xy_1, \ldots, x^{n_1-1}y_1$ respectively, and so on. Of course we assume a certain ordering $\mu_1, \ldots, \mu_k$ of different eigenvalues is fixed.

\begin{definition}
Given a non-constant $x\in K$, let $v_1, \ldots, v_n$ be a canonical fundamental system and let $y \in K$ be such that $\Delta(x,y)=0$. Then $y$ (or the pair $(x,y)$) is said to be a \emph{proper} solution if $y=\sum a_i v_i$ with $a_i \in C^{\times}$, that is, if $y$ is not in the linear span of a proper subset of $\{v_1, \ldots, v_n\}$.
\end{definition}

A solution is proper if and only if it does not satisfy a linear differential equation of lower order.

\begin{lemma}\label{proper-criterion}
A pair $(x,y) \in K^2$ is a proper solution to \eqref{main-eq} if and only if $y, \partial_x y, \ldots, \partial_x^{n-1}y$ are $C$-linearly independent.
\end{lemma}
\begin{proof}
Let $v_1, \ldots, v_n$ be as above and $y=\sum a_i v_i$. Since $v_1, \ldots, v_n$ are $C$-linearly independent, the Wronskian $W(\bar{v})=\det(\partial_x^l v_i)$ is non-zero. It is easy to check that $\partial_x^l(v_i) = f_{li}(x) v_i$ where $f_{li}$ is a rational function over $\mathbb{Q}(\mu_1, \ldots, \mu_k)$.  Furthermore, none of the $f_{li}(x)$ is zero (as $x$ is non-constant). Let $H_x$ be the $n \times n$ matrix with entries $f_{li}(x)$. Then $W(\bar{v}) = \det(H_x) \cdot \prod_{i=1}^m v_i$. 
Consider the following system of equations with respect to $v$'s:
$$\partial_x^l (y) = \sum_{i=1}^m a_i f_{li}(x) v_i,~ l=0, \ldots, n-1.$$
Its determinant is $\det(H_x) \cdot \prod_{i=1}^m a_i$
which is non-zero if and only if none of the $a_i$'s is zero. This finishes the proof.
\end{proof}

Let $(x,y)$ be a proper solution. Then we can assume $y=v_1+\ldots+v_n$. Let $H_x$ be as in the proof and denote its rows by $H_x^l$. It is an invertible linear transformation of $K^n$. Let $L_x$ be its inverse with coordinate functions (rows) $L_x^i:K^n \rightarrow K$. Thus $$\partial_x^l(y) = H_x^l(v_1,\ldots, v_n) \mbox{ and } v_i = L_x^i(y, \partial_xy, \ldots, \partial_x^{n-1}y).$$
It is also worth mentioning that when $p(\lambda)$ does not have multiple roots, $H_x$ and $L_x$ do not depend on $x$, they depend only on $\lambda_i$'s. 
Note also that if $\Delta(x,y)=0$ and $x$ is non-constant then $\Delta(x, \partial_x y)=0$. In particular, if $(x,y)$ is a proper solution then $y, \partial_xy, \ldots, \partial_x^{n-1}y$ form a fundamental system of solutions.
These considerations will be useful in Section \ref{Rotundity and freeness}.

Now we are ready to prove the Ax-Schanuel inequality for \eqref{main-eq}.
\begin{theorem}\label{Ax-Schanuel-for-higher-order}
Let $(x_i,y_i),~ i=1,\ldots,m,$ be proper solutions to the equation \eqref{main-eq} in $\mathcal{K}$. Then
\begin{equation}\label{inequality-for-higher-order}
\td_CC(\bar{x},\bar{y},\partial_{\bar{x}}\bar{y},\ldots,\partial_{\bar{x}}^{n-1}\bar{y}) \geq \ldim_{\mathbb{Q}}(\lambda_1 \bar{x},\ldots,\lambda_n \bar{x}/C)+1,
\end{equation}
where $\partial^j_{\bar{x}}\bar{y}=(\partial^j_{x_1}y_1, \ldots, \partial^j_{x_m}y_m)$
\end{theorem}

In particular, if we assume $\lambda_1 \bar{x},\ldots,\lambda_n \bar{x}$ are $\mathbb{Q}$-linearly independent modulo $C$ then $\td_CC(\bar{x},\bar{y},\partial_{\bar{x}}\bar{y},\ldots,\partial_{\bar{x}}^{n-1}\bar{y})\geq mn+1$. This is possible only if $\lambda_1,\ldots,\lambda_n$ are linearly independent over $\mathbb{Q}$. In fact we can always assume it is the case; otherwise both the transcendence degree and the linear dimension will decrease and we will be reduced to the same inequality for a smaller $n$.
Note also that the case $n=1$ is exactly Ax's theorem for the exponential differential equation.

\begin{proof}[Proof of Theorem \ref{Ax-Schanuel-for-higher-order}]
For each $i$ let $v_{ij} \in K^{\times},~ j=1,\ldots,n,$ be a canonical fundamental system of solutions to $\Delta(x_i, y)=0$. Then for every $i$ the $C$-linear span of $v_{i1},\ldots,v_{in}$ is the same as that of $y_i,\partial_{x_i}y_i,\ldots,\partial_{x_i}^{n-1}y_i$, for $(x_i,y_i)$ is a proper solution. In particular, the field $C(\bar{x},\bar{y},\partial_{\bar{x}}\bar{y},\ldots,\partial_{\bar{x}}^{n-1}\bar{y})$ is equal to the field extension of $C$ generated by $\bar{x}$ and all the $v_{ij}$. Therefore
\begin{align*}
\td_CC(\bar{x},\bar{y},\partial_{\bar{x}}\bar{y},\ldots,\partial_{\bar{x}}^{n-1}\bar{y}) & = \td_CC(\mu_1 \bar{x},\ldots,\mu_k \bar{x},\bar{v}_1,\ldots,\bar{v}_n)\\
& \geq \ldim_{\mathbb{Q}}(\mu_1 \bar{x},\ldots,\mu_k \bar{x}/C)+1 \\
& = \ldim_{\mathbb{Q}}(\lambda_1 \bar{x},\ldots,\lambda_n \bar{x}/C)+1
\end{align*}
where $\bar{v}_j$ is the tuple $(v_{1j},v_{2j},\ldots,v_{mj})$. The inequality follows from Ax's theorem applied to the tuple $(\mu_1 \bar{x},\ldots,\mu_k \bar{x})$ taking into account that $(\bar{v}_1,\ldots,\bar{v}_n)$ contains a solution $y_{ij}$ for each of the equations $\Exp(\mu_i x_j, y_{ij})$. 
\end{proof}

{Let us also note that one can prove (with a similar argument) an analogue of Theorem \ref{Ax-Schanuel-for-higher-order} for fields with several commuting derivations (using the corresponding version of Ax's theorem). Nevertheless, we prefer working in differential fields with a single derivation for simplicity and do not consider a multi-derivative version of the above theorem.}




\section{The complete theory}\label{Existential-closedness}
\setcounter{equation}{0}

Having established a predimension inequality (see Section \ref{Predimension} for details on predimensions) for higher order linear differential equations, we want to find an appropriate existential closedness property and thus give an axiomatisation of the complete theory of the corresponding reducts. 

First, let us see which language we should work in. An obvious option would be simply taking a binary predicate for the solutions of the equation \eqref{main-eq}. But the inequality \eqref{inequality-for-higher-order} cannot be written as a first order statement (axiom scheme) in this language. This is because derivatives of $y_i$'s are involved in \eqref{inequality-for-higher-order}. Therefore we need to take a predicate of higher arity which will have variables for the derivatives of $y$'s as well. Thus we will work in the language $\mathfrak{L}_{\E_n} = \{ +,\cdot, \E_n, 0, 1, \lambda_1, \ldots, \lambda_n \} $ where $\lambda_1, \ldots, \lambda_n$ are constant symbols for the eigenvalues and $\E_n(x, z_0, z_1,\ldots,z_{n-1})$ is an $(n+1)$-ary predicate. It will be interpreted in a differential field $\mathcal{K}$ as the set $$\left\{ (x,\bar{z}) \in K^{n+1} : \exists z_n \left[z_n+\sum_{i=0}^{n-1} c_iz_i=0 \wedge \bigwedge_{i=0}^{n-1}\D z_i = z_{i+1} \D x\right] \right\}.$$ Note that since $\lambda_1, \ldots, \lambda_n$ are in the language, the coefficients $c_0, \ldots, c_{n-1}$ are $\emptyset$-definable. The theory of reducts of differentially closed fields to the language $\mathfrak{L}_{\E_n}$ will be denoted by $T_{\E_n}$. Also the field of constants can be defined as {$C=\{ c : \E_n(c,0,1,0, \ldots, 0) \}$.} 

Observe that $\Exp$ can be defined in an $\E_n$-reduct of a differential field. {Namely,
$$\mathcal{K} \models \Exp(\lambda_i x, y) \leftrightarrow \E_n(x,y, \lambda_iy, \ldots, \lambda_i^{n-1}y)$$
for any $i \in \{1, \ldots, n\}$. Indeed, if $\Exp(\lambda_i x, y)$ holds then by Lemma \ref{solutions} $\Delta(x,y)=0$ and $\partial^j_x y = \lambda_i^j y$ for each $j$ and so $(x,y, \lambda_iy, \ldots, \lambda_i^{n-1}y) \in \E_n$. Conversely, if  $(x,y, \lambda_iy, \ldots, \lambda_i^{n-1}y) \in \E_n$ then $y$ can be written as a $C$-linear combination of the fundamental system of solutions. Moreover, we must have $\partial^j_x y = \lambda_i^j y$ for $j=0,\ldots,n-1$. This system of equations implies that $\partial_x y = \lambda_i y$ and so $\Exp(\lambda_i x, y)$ holds.}

In fact $\Exp$ and $\E_n$ are interdefinable. So we can just translate the axiomatisation for the exponential differential equation to the language $\mathfrak{L}_{\E_n}$ and get an axiomatisation of $T_{\E_n}$. However we want an axiomatisation based on the Ax-Schanuel inequality proved in Section \ref{Higher-order-section}. In other words, we want to understand which systems of equations in $\mathfrak{L}_{\E_n}$-reducts of differentially closed fields have solutions, and prove that \eqref{inequality-for-higher-order} is an adequate predimension inequality.

\begin{notation}
If $\partial_xy_i = \mu_i y_i$ then let $g_{ijl}(X)$ be the algebraic polynomial (over $\mathbb{Q}(\mu_i)$) for which $\partial_x^l(x^j y_i)=g_{ijl}(x)y_i$. In particular $g_{i0l}= \mu_i^l$. Also denote $N_i:=1+\sum_{j<i} n_j$.
\end{notation}

Now we formulate a number of axioms and axiom schemes for an $\mathfrak{L}_{\E_n}$-structure $\mathcal{F}_{\E_n}$. In such a structure we let $\Exp(x,y)$ be the relation defined by the formula $\E_n(\lambda_1^{-1}x,y, \lambda_1y, \ldots, \lambda_1^{n-1}y)$.

\begin{itemize}
\item[A1'] $F$ is an algebraically closed field. 

\item[A2']$C:= \{ c\in F: \mathcal{F}_{\E_n}\models \E_n(c,1, 0, \ldots, 0) \}$ is an algebraically closed subfield of $F$ and $\lambda_1,\ldots, \lambda_n$ are non-zero elements of $C$ satisfying the appropriate algebraic relations. In particular $\lambda_{N_i}=\lambda_{N_i+1}=\ldots=\lambda_{N_i+n_i-1} =: \mu_i$ for every $i$. 

\item[A3'] $\E_n(x, z_0, \ldots, z_{n-1})$ holds if and only if 
there are $y_1, \ldots, y_k \in F^{\times}$ with $\Exp(\mu_i x, y_i)$ and elements $a_{ij} \in C$ such that
$$z_l = \sum_{i=1}^k \sum_{j=0}^{n_i-1} a_{ij}g_{ijl}(x)y_i,$$
for $l=0, \ldots, n-1$.


\item[A4'] $\Exp(F) = \{ (x,y)\in F^2: \Exp(x,y)\}$ is a subgroup of $G_1(F)$ containing $G_1(C)$.

\item[A5'] The fibres of $\Exp$ in $\mathbb{G}_a(F)$ and $\mathbb{G}_m(F)$ are cosets of the subgroups $\mathbb{G}_a(C)$ and $\mathbb{G}_m(C)$ respectively.

\item[AS'] Let $x_i,z_{ij}\in F \setminus C,~ 1\leq i \leq m, 0\leq j <n,$ be such that $z_{i0}, \ldots, z_{i,n-1}$ are {$C$-linearly} independent and
$$\mathcal{F}_{\E_n} \models \bigwedge_i \E_n(x_i, z_{i0}, \ldots, z_{i,n-1}).$$ Then for each $1 \leq d \leq mn$ if $\ldim_{\mathbb{Q}}(\lambda_1 \bar{x},\ldots,\lambda_n \bar{x}/C) \geq d$ then
$$\td_CC(\bar{x},\bar{z}_0,\bar{z}_1,\ldots,\bar{z}_{n-1}) \geq d+1.$$

\item[{NT'}] {$F \supsetneq C$}.
\end{itemize}

{As in the case of AS, a compactness argument can be used here as well to show that AS' can be expressed as a first-order axiom scheme. Observe also that $C$-linear independence of $z_{i0}, \ldots, z_{i,n-1}$ is first-order here (as opposed to $\mathbb{Q}$-linear independence), for $C$ is definable in the reduct.}

\begin{definition}
{An $\E_n$-field is a model of A1'-A5',AS'.} 
\end{definition}


\begin{lemma}\label{Axn-implies-Ax1}
Let $\mathcal{F}_{\E_n}$ be a model of \emph{A1'-A5'}. Then it satisfies \emph{AS'} iff the relation $\Exp(x,y)$ satisfies \emph{AS}.
\end{lemma}
\begin{proof}
Let $x_1,\ldots, x_m \in F$ be $\mathbb{Q}$-linearly independent modulo $C$. Then $\mu_1 x_1,\ldots, \mu_1 x_m$ are such as well. Denote $\mu_s x_i =: u_{m(s-1)+i}$ for $i=1, \ldots, m,~ s=1, \ldots, k$. If $\ldim_{\mathbb{Q}}(u_1,\ldots, u_{mk}/C) = d \geq m$ then assume without loss of generality that $u_1,\ldots, u_d$ are linearly independent over the rationals modulo $C$. Let $v_i \in F$ be such that $\mathcal{F}_{\E_n} \models \Exp(u_i, v_i)$. Then AS' implies that 
$$\td_CC(x_1,\ldots, x_m, v_1, \ldots, v_{mk}) \geq d+1.$$
For each $i>d$ there are integers $m_i, m_{i1},\ldots, m_{id}$ such that $m_i u_i + m_{i1}u_1+\ldots + m_{id}u_d = c \in C$.

Denote $v= v_i^{m_i} \prod_{j=1}^d v_j^{m_{ij}}$. By A4' we have $\Exp(c,v)$. But also $\Exp(c,1)$ holds and using A5' we deduce that $v \in C$. Hence $v_1, \ldots, v_d, v_i$ are algebraically dependent over $C$. Therefore 
$$\td_CC(\bar{x}, v_1,\ldots, v_d) \geq d+1.$$
Now we can easily deduce that $\td_CC(\bar{x}, v_1,\ldots, v_m) \geq m+1$ and we are done.

The converse follows from the proof of Theorem \ref{Ax-Schanuel-for-higher-order}.
\end{proof}


\begin{notation}
Let $\pr_j : K^{m(n+1)} \rightarrow K^{2m}$ be defined as
$$\pr_j:(\bar{x},\bar{v}_0, \ldots, \bar{v}_{n-1}) \mapsto (\bar{x},\bar{v}_j),$$ where $\bar{v}_j = (v_{1j}, \ldots, v_{mj})$.

Also we will denote the set $\{ (\bar{x},\bar{z}_0,\ldots,\bar{z}_{n-1}) \in F^{m(n+1)} : \mathcal{F}_{\E_n} \models \bigwedge_{i=1}^m \E_n(x_i,\bar{z}^i) \}$ by $\E^m_n(F)$ where $\bar{z}^i=(z_{i0},\ldots,z_{i,n-1})$.
\end{notation}

\begin{definition}\label{rotundity-def}
An irreducible variety $ V \subseteq K^{m(n+1)}$ is called $\E_n$-$\Exp$-\textit{rotund} if $V_1:=\pr_1 (V) \subseteq G_{m}(K)$ is $\Exp$-rotund and 
\begin{equation}\label{rotund-eq}
(\bar{x}, \bar{y}) \in V_1 \Longrightarrow
(\bar{x}, \bar{y}, \mu_1 \bar{y}, \ldots, \mu_1^{n-1} \bar{y}) \in V.
\end{equation}
\end{definition}
We could of course replace $\mu_1$ in \eqref{rotund-eq} by any $\mu_i$. As $\Exp$-rotundity is a definable property, so is $\E_n$-$\Exp$-rotundity (in the language of rings). 

Now we formulate the existential closedness property for an $\E_n$-field $\mathcal{F}_{\E_n}$. 

\begin{itemize}
\item[EC'] For each irreducible $\E_n$-$\Exp$-rotund variety $ V \subseteq F^{m(n+1)}$ the intersection $V(F) \cap \E^m_n(F)$ is non-empty.
\end{itemize}


This statement can be given by a first-order axiom scheme in the language $\mathfrak{L}_{\E_n}$, for $\E_n$-$\Exp$-rotundity is a first-order property. 

\begin{lemma}\label{EC1-holds-DCF}
If $\mathcal{K}$ is a differentially closed field then $\mathcal{K}_{\E_n}$ satisfies \emph{EC'}.
\end{lemma}
\begin{proof}
Let $ V \subseteq K^{m(n+1)}$ be an $\E_n$-$\Exp$-rotund variety. Then $V_1= \pr_1(V)$ is an $\Exp$-rotund variety. So by Theorem \ref{exp-axioms} and Lemma \ref{exp-lemma} there is a point $(\bar{x}, \bar{y}) \in V_1$ such that $\mathcal{K}_{\E_n} \models \Exp(\mu_1 x_i, y_i)$ for each $i=1,\ldots, m$. By definition we have
$$(\bar{x}, \bar{y}, \mu_1 \bar{y}, \ldots, \mu_1^{n-1} \bar{y}) \in V.$$
It is also clear that $$\mathcal{K}_{\E_n} \models \E_n(\bar{x}, \bar{y}, \mu_1 \bar{y}, \ldots, \mu_1^{n-1} \bar{y})$$ and we are done.
\end{proof}



\begin{lemma}\label{higher-order-EC}
If $\mathcal{F}_{\E_n}$ satisfies \emph{A1'-A5', AS', EC'} then $\Exp(x,y)$ satisfies \emph{EC}.
\end{lemma}
\begin{proof}
Suppose $W \subseteq G_m(F)$ is an $\Exp$-rotund variety defined over a set $A \subseteq F$. Let $\mathbb{F} \supseteq F$ be a saturated algebraically closed field and pick $(\bar{a}, \bar{b}) \in \mathbb{F}^{2n}$ a generic point of $W$. Let $V \subseteq F^{m(n+1)}$ be the algebraic locus of $(\mu_1^{-1}\bar{a}, \bar{b}, \mu_1 \bar{b}, \ldots, \mu_1^{n-1}\bar{b})$ over $A\mu_1$. Then $V$ is $\E_n$-$\Exp$-rotund and hence $ V(F)\cap \E_n^m(F) \neq \emptyset$. By our construction of $V$ we also know that any point in that intersection must be of the form $(\mu_1^{-1}\bar{x}, \bar{y}, \mu_1 \bar{y}, \ldots, \mu_1^{n-1} \bar{y})$. Then $(\bar{x}, \bar{y}) \in W$ and by A3' $\mathcal{F}_{\E_n} \models \Exp(\bar{x}, \bar{y})$. So $W(F) \cap \Exp^n(F) \neq \emptyset$.
\end{proof}

Finally, we can deduce that the given axioms form a complete theory.

\begin{theorem}\label{higher-order-complete-axioms}
The axioms and axiom schemes \emph{A1'-A5', AS', EC', NT'} axiomatise the complete theory $T_{\E_n}$.
\end{theorem}
\begin{proof}
Indeed, Lemmas \ref{Axn-implies-Ax1}, \ref{EC1-holds-DCF} and \ref{higher-order-EC} show that an $\mathfrak{L}_{\E_n}$-structure $\mathcal{F}_{\E_n}$ satisfies A1'-A5', AS', and EC' if and only if the relation $\Exp(x,y)$ satisfies the axioms A1-A4, AS, and EC. The latter collection of axioms axiomatises the theory $T_{\Exp}$ by Theorem \ref{exp-axioms}. Now the desired result follows as the relations $\Exp$ and $\E_n$ are interdefinable due to A3'.
\end{proof}

\section{Rotundity and freeness}\label{Rotundity and freeness}
\setcounter{equation}{0}

Though EC' is an appropriate existential closedness property for $\E_n$-fields, our definition of $\E_n$-$\Exp$-rotundity is not that natural. Indeed, the inequality given by AS' is not reflected in it and also the notion of $\E_n$-$\Exp$-rotundity is far from being a necessary condition for a variety to intersect $\E_n$. As we saw, $\E_n$-$\Exp$-rotund varieties have a very special $\E_n$-point, which is essentially (constructed from) an exponential point. For these reasons we define another notion of rotundity (and strong rotundity), {called $\E_n$-rotundity (respectively, strong $\E_n$-rotundity)}\footnote{{Dropping the $\Exp$ from the term $\E_n$-$\Exp$-rotundity indicates that the new notion of rotundity is ``more independent'' from $\Exp$.}}, which will be more intuitive and natural. (That definition will not be as simple as Definition \ref{rotundity-def} though.) We will  see in particular that strongly $\E_n$-rotund varieties will contain proper $\E_n$-points.

Recall that $\mu_1, \ldots, \mu_k$ are the distinct eigenvalues of our differential equation. As before, for $z_{ij}, i=1, \ldots, m, j=0, \ldots, n-1,$ denote $\bar{z}^i=(z_{i0},\ldots,z_{i,n-1})$ and $\bar{z}_j=(z_{1j},\ldots,z_{mj})$.
Define the map $$\tilde{L}: K^{m(n+1)} \rightarrow K^{m(n+1)}$$ by
$$\tilde{L}: (\bar{x},\bar{z}_0, \ldots,\bar{z}_{n-1}) \mapsto (\bar{x},L^1_{x_1}(\bar{z}^1), \ldots, L^1_{x_m}(\bar{z}^m), \ldots, L^n_{x_1}(\bar{z}^1), \ldots, L^n_{x_m}(\bar{z}^m)),$$
where $L^j_{x_i}$ is as in Section \ref{Higher-order-section}. 
Let $\tilde{H}$ be its inverse map. 
Recall that for $1 \leq i \leq k$ we denoted $N_i=1+ \sum_{j<i} n_j$.
Define maps $R:F^{m(n+1)} \rightarrow F^{m(k+1)}$ and $\tilde{R}:F^{m(n+1)} \rightarrow F^{2km}$ as follows:
$$R: (\bar{x},\bar{v}_1, \ldots,\bar{v}_{n}) \mapsto (\bar{x}, \bar{v}_{N_1}, \ldots,\bar{v}_{N_k}),$$
$$\tilde{R} :(\bar{x},\bar{v}_1, \ldots,\bar{v}_{n}) \mapsto (\mu_1\bar{x}, \ldots, \mu_k\bar{x}, \bar{v}_{N_1}, \ldots,\bar{v}_{N_k}).$$

\begin{definition}
An irreducible variety $ V \subseteq F^{m(n+1)}$ is called (\emph{strongly}) $\E_n$-\textit{rotund} if $V':= \tilde{R}\circ\tilde{L} (V) \subseteq G_{km}(F)$ is (strongly) $\Exp$-rotund and $V'':=R(\tilde{L}(V))$ satisfies the following property: 
\begin{equation*}
(\bar{x}, \bar{y}_1, \ldots, \bar{y}_k) \in V'' \Rightarrow
\tilde{H}(\bar{x}, \bar{y}_1, x\bar{y}_1, \ldots, x^{n_1-1}\bar{y}_1, \ldots, \bar{y}_k, x\bar{y}_k, \ldots, x^{n_k-1}\bar{y}_k) \in V.
\end{equation*}
\end{definition}

One can use this notion of rotundity to formulate an appropriate existential closedness statement (that is, the above EC' but for $\E_n$-rotund varieties instead of $\E_n$-$\Exp$-rotund ones) which, with A1'-A5' and AS', axiomatises $T_{\E_n}$. The following result shows that this notion of rotundity fits better with our differential equation.

\begin{proposition}
Let $\mathcal{K}$ be a differentially closed field. If $V \subseteq K^{m(n+1)}$ is a strongly $\E_n$-rotund variety defined over $C$ then $V(K)$ has a proper $\E_n$-point.
\end{proposition}
\begin{proof}
Indeed strong $\Exp$-rotundity of $V'$ implies that it has a non-constant $\Exp$-point. This point obviously gives rise to a proper $\E_n$-point on $V$.
\end{proof}


In sufficiently saturated models of $T_{\E_n}$ every (strongly) $\E_n$-rotund variety contains a generic (proper) $\E_n$-point. The converse holds for ``free'' varieties: (strong) $\E_n$-rotundity is a necessary condition for a free variety to have a generic (proper) $\E_n$-point. We give precise definitions below. 

{In the rest of the section we assume that these are linearly independent
over $\mathbb{Q}$. The general case involves no significant additional ideas (after
taking a basis), but is notationally messier.}


\begin{definition}\footnote{Note that if we do not require $\mu_1, \ldots, \mu_k$ to be linearly independent over $\mathbb{Q}$ (as we assumed) then the above definition would not make sense. Of course in that case we could just change the definition of the map $R$ appropriately and get the same notion of freeness.}
An irreducible variety $ V \subseteq F^{m(n+1)}$ is called $\E_n$-\textit{free} if $V':= \tilde{R}\circ\tilde{L} (V) \subseteq G_{km}(F)$ is $\Exp$-free.
\end{definition}

The next result follows from Proposition \ref{exp-free} and some standard observations on generic points. It can also be proven using Theorem \ref{Ax-Schanuel-for-higher-order}.

\begin{proposition}
Suppose $V \subseteq F^{m(n+1)}$ is an irreducible and free variety defined over $C$. If $V$ has a proper generic (over $C$) $\E_n$-point then it must be strongly $\E_n$-rotund.
\end{proposition}

\section{Predimension}\label{Predimension}
\setcounter{equation}{0}

We denote by $T_{\Exp}^0$ the $\mathfrak{L}_{\Exp}$-theory given by the axioms A1-A4, AS. Similarly, $T_{\E_n}^0$ is the $\mathfrak{L}_{\E_n}$-theory consisting of the axioms A1'-A5', AS'. The results of Section \ref{Existential-closedness} show that $T^0_{\Exp}$ and $T^0_{\E_n}$ (as well as $T_{\Exp}$ and $T_{\E_n}$)  are essentially the same theory given in two different languages. In particular, every model $\mathcal{F}_{\E_n}$ of $T^0_{\E_n}$ (or $T_{\E_n}$) can be canonically made into a model $\mathcal{F}_{\Exp}$ of $T^0_{\Exp}$ (respectively $T_{\Exp}$) and vice versa. This relationship allows us to deduce that the predimension inequality \eqref{inequality-for-higher-order} is adequate. We proceed towards this goal in this section. 

We will first prove that an embedding of $\E_n$-fields is the same as an embedding of the corresponding $\Exp$-fields. 
For an $\E_n$-field $\mathcal{F}_{\E_n}$ (or an $\Exp$-field $\mathcal{F}_{\Exp}$) we let $C_F$ denote its field of constants.

\begin{lemma}\label{embeddings}
{Let $\mathcal{K}_{\E_n}$ and $\mathcal{F}_{\E_n}$ be two $\E_n$-fields with an embedding of fields $f: K \hookrightarrow F$. Then $f: \mathcal{K}_{\E_n} \hookrightarrow \mathcal{F}_{\E_n}$ is an embedding of $E_n$-fields if and only if it is an embedding of $\Exp$-fields $f:\mathcal{K}_{\Exp} \hookrightarrow \mathcal{F}_{\Exp}$.}
\end{lemma}
\begin{proof}
Since $\Exp$ is quantifier-free definable in an $\E_n$-field, we only need to show that an embedding of $\Exp$-fields is also an embedding of the corresponding $\E_n$-fields. Identifying $K$ with $f(K)$ we can assume $\mathcal{K}_{\Exp}\subseteq \mathcal{F}_{\Exp}$. Let $a, b_0, \ldots, b_{n-1} \in K$ be such that
$$\mathcal{F}_{\Exp} \models \E_n(a, b_0, \ldots, b_{n-1}).$$ We shall show that $$\mathcal{K}_{\Exp} \models \E_n(a, b_0, \ldots, b_{n-1}).$$

We can assume that $a$ is non-constant. By A3' we know that there are $e_1, \ldots, e_k \in F^{\times}$ with $\Exp(\mu_i a, e_i)$ and elements $a_{ij} \in C_F$ such that
$$b_l = \sum_{i=1}^k \sum_{j=0}^{n_i-1} a_{ij}g_{ijl}(a)e_i,$$
for $l=0, \ldots, n-1$ (here $g_{ijl}$ is as in Section \ref{Existential-closedness}). If $\Exp(u,v)$ holds for some $u,v$ then $\Exp(u, cv)$ holds as well for any constant $c$. Hence we can assume without loss of generality that $a_{ij}$ is either $0$ or $1$. As $g_{ijl}(X) \in C_K[X]$, we can express all $e_i$'s with $a_{i0} =1$ in terms of $g_{ijl}(a)$ and $b_l$ (this is because the corresponding determinant does not vanish). Hence $e_i \in K$ and we are done by A3' again.
\end{proof}

This lemma shows that the category of $\E_n$-fields with morphisms being embeddings is isomorphic to the category of $\Exp$-fields again with embeddings as morphisms.

Let $C$ be a countable algebraically closed field with $\td(C/\mathbb{Q})=\aleph_0$ and let $\mathfrak{C}^{\Exp}$ be the class of all countable models of $T^0_{\Exp}$ having $C$ as field of constants. Further, denote the subclass of $\mathfrak{C}^{\Exp}$, consisting of $\Exp$-fields which have finite transcendence degree over $C$ as a field, by $\mathfrak{C}^{\Exp}_{f.g.}$.\footnote{The subscript $f.g.$ here stands for \textit{finitely generated}, and it means that the structures of $\mathfrak{C}^{\Exp}_{f.g}$ are finitely generated as algebraically closed fields over $C$, that is, they have finite transcendence degree over $C$.} 

Following \cite{Kirby-semiab} for $\mathcal{K}_{\Exp} \in \mathfrak{C}^{\Exp}_{f.g.}$ (with domain $K$) define 
\begin{align*}
\sigma_{\Exp}(\mathcal{K}_{\Exp}):=\max \{ n: &\mbox{ there are } a_i, b_i \in K,~ i=1,\ldots, n, \mbox{ with } a_i\mbox{'s}\\ & \mbox{ linearly independent over } \mathbb{Q} \mbox{ mod } C \mbox{ and } \mathcal{K}_{\Exp} \models \Exp(a_i,b_i) \}
\end{align*}
and
$$\delta_{\Exp}(\mathcal{K}_{\Exp}):=\td_{C}(\mathcal{K}_{\Exp}) - \sigma_{\Exp}(\mathcal{K}_{\Exp}).$$

Firstly note that $\sigma_{\Exp}$ is well defined and finite since the Ax-Schanuel inequality bounds the number $n$ in consideration by $\td_{C}C(\bar{a},\bar{b})$ which, in its turn, is bounded by $\td_{C}(\mathcal{K}_{\Exp})$.

Secondly, the Ax-Schanuel inequality is equivalent to saying that $\delta_{\Exp}(\mathcal{K}_{\Exp}) \geq 0$ for all $\mathcal{K}_{\Exp} \in \mathfrak{C}^{\Exp}_{f.g.}$ where equality holds if and only if $K=C$.

The function $\delta_{\Exp}$ is a predimension in the sense of Hrushovski (see \cite{Kirby-semiab, Aslanyan-thesis}). In particular, it can be used to define strong extensions.

\begin{definition}\label{strong-exp}
Let $A, B \in \mathfrak{C}^{\Exp}$ with $A\subseteq B$. Then we say $A$ is a \emph{strong substructure} of $B$ or $B$ is a \emph{strong extension} of $A$, denoted $A \leq B$, if for every $X \in \mathfrak{C}^{\Exp}_{f.g.}$ with $X \subseteq B$ we have $\delta_{\Exp}(X \cap A) \leq \delta_{\Exp}(X)$.
\end{definition}

Kirby proved in \cite{Kirby-semiab} that the class $\mathfrak{C}^{\Exp}$ satisfies the \textit{strong amalgamation property} and is an $\aleph_0$-\emph{amalgamation category} (see \cite{Droste-Gobel,Kirby-thesis,Kirby-semiab} for details). So one can carry out an (uncollapsed) Hrushovski construction and end up with a unique (up to isomorphism) countable $\Exp$-field $\mathcal{U}_{\Exp}$ which is universal and saturated with respect to strong embeddings. 
This \textit{strong amalgam} $\mathcal{U}_{\Exp}$ also satisfies the \emph{strong existential closedness} property, that is, it is existentially closed in strong extensions. Furthermore, the following result holds showing that the Ax-Schanuel inequality is adequate for $\Exp$. 

\begin{theorem}[\cite{Kirby-thesis,Kirby-semiab}]\label{Adequacy-for-exp}
Let $\mathcal{F}$ be the countable saturated differentially closed field. Then $\mathcal{U}_{\Exp}$ is isomorphic to $\mathcal{F}_{\Exp}$. In particular, the $\Exp$-reduct of any differentially closed field is elementarily equivalent to $\mathcal{U}_{\Exp}$.
\end{theorem}

In order to define the predimension for $\E_n$-fields we first observe that Theorem \ref{Ax-Schanuel-for-higher-order} can be reformulated to give a lower bound for transcendence degree not only for proper solutions but for arbitrary ones. Recall that $\mu_1,\ldots,\mu_k$ are all the distinct eigenvalues of our equation with multiplicities $n_1, \ldots, n_k$ respectively. Let $v_1,\ldots,v_n$ be a canonical fundamental system of solutions of $\Delta(x,y)=0$. 
For a solution $(x,y)$ (with $x$ non-constant) we have a unique representation $y=c_1v_1+\ldots + c_n v_n$ with $c_i \in C$. 
For $1\leq i \leq k$  we define 
$$\epsilon_i(y) := 
\begin{cases}
1, \text{ if for some } j \text{ with } N_i\leq j < N_{i+1} \text{ we have } c_j \neq 0, \\
0, \text{ otherwise}.
\end{cases}
$$
Then set $\epsilon(y):=(\epsilon_1(y),\ldots,\epsilon_k(y))$ and denote $\epsilon(y)x:=(\epsilon_1(y)\mu_1x,\ldots,\epsilon_k(y)\mu_kx)$.

Now it is easy to see that Theorem \ref{Ax-Schanuel-for-higher-order} is equivalent to the following.

\begin{theorem}\label{Ax-Schanuel-for-higher-order-2}
Let $(x_i,y_i),~ i=1,\ldots,m,$ be solutions to the equation \eqref{main-eq} in $\mathcal{K}$ with $x_i \in K \setminus C$. Then
\begin{equation}\label{inequality-for-higher-order-2}
\td_CC(\bar{x},\bar{y},\partial_{\bar{x}}\bar{y},\ldots,\partial_{\bar{x}}^{n-1}\bar{y}) - \ldim_{\mathbb{Q}}(\epsilon(y_1)x_1, \ldots, \epsilon(y_m)x_m/C)\geq 1.
\end{equation}
\end{theorem}


Now we define the predimension. As above, fix a countable algebraically closed field $C$ with $\td(C/\mathbb{Q})=\aleph_0$ and let $\mathfrak{C}^{\E_n}$ be the collection of all $\E_n$-fields with field of constants $C$. Further, $\mathfrak{C}^{\E_n}_{f.g.}$ consists of those structures from $\mathfrak{C}^{\E_n}$ which have finite transcendence degree over $C$. For $\mathcal{K}_{\E_n} \in \mathfrak{C}^{\E_n}_{f.g.}$ (with domain $K$) define 
\begin{align*}
\sigma_{\E_n}(\mathcal{K}_{\E_n}):=\max \{ & \ldim_{\mathbb{Q}}(\epsilon^1a_1, \ldots, \epsilon^ma_m/C):  \mbox{ where } \epsilon^i \in \{0, 1\}^k,~ a_i \in K  \mbox{ such that}\\ 
&\mbox{there are } b^0_i, \ldots, b^{n-1}_i \in K \mbox{ with } \epsilon(b^0_i)=\epsilon^i \mbox{ and} \\
 &\mathcal{K}_{\E_n} \models \E_n(a_i,b^0_i,\ldots,b^{n-1}_i),~ i=1,\ldots, m \}
\end{align*}
and
$$\delta_{\E_n}(\mathcal{K}_{\E_n}):=\td_{C}(\mathcal{K}_{\E_n}) - \sigma_{\E_n}(\mathcal{K}_{\E_n}).$$

Then the inequality \eqref{inequality-for-higher-order-2} states precisely that $\delta_{\E_n}(\mathcal{K}_{\E_n})\geq 0$ for all $\mathcal{K}_{\E_n} \in \mathfrak{C}_{f.g.}$ and equality holds if and only if $\mathcal{K}_{\E_n}=C$. 
One can also define strong embeddings of $\E_n$-fields by replacing $\delta_{\Exp}$ with $\delta_{\E_n}$ in Definition \ref{strong-exp}.

\begin{lemma}
For an $\E_n$-field $\mathcal{K}_{\E_n} \in \mathfrak{C}_{f.g.}$ we have
$$\sigma_{\E_n}(\mathcal{K}_{\E_n})=\sigma_{\Exp}(\mathcal{K}_{\Exp}),~ \delta_{\E_n}(\mathcal{K}_{\E_n})=\delta_{\Exp}(\mathcal{K}_{\Exp}).$$
Hence an embedding of $\E_n$-fields $\mathcal{K}_{\E_n} \hookrightarrow \mathcal{F}_{\E_n}$ is strong if and only if it is strong as an embedding of the corresponding $\Exp$-fields $\mathcal{K}_{\Exp} \hookrightarrow \mathcal{F}_{\Exp}$. Furthermore, the category of $\E_n$-fields with morphisms being strong embeddings is isomorphic to the category of $\Exp$-fields with strong embeddings. 
\end{lemma}
\begin{proof}
We only need to show that $\sigma_{\E_n}(\mathcal{K}_{\E_n})=\sigma_{\Exp}(\mathcal{K}_{\Exp})$. Let $m=\sigma_{\Exp}(\mathcal{K_{\Exp}})$ and $(a_1,\ldots,a_m)\in K^m$ be linearly independent over $\mathbb{Q}$ mod $C$ such that for some $b_1,\ldots,b_m \in K^{\times}$ we have $\mathcal{K}_{\Exp} \models \Exp(\mu_1 a_i,b_i)$ for each $i$. Then $\mathcal{K}_{\E_n} \models \E_n(a_i, b_i, \mu_1 b_i, \ldots, \mu_1^{n-1} b_i)$. Clearly, $\epsilon(b_i)=1$ only for $i=1$. Hence $\ldim_{\mathbb{Q}}(\epsilon^1a_1, \ldots, \epsilon^ma_m/C)=m$ and so $\sigma_{\E_n}(\mathcal{K}_{\E_n}) \geq m=\sigma_{\Exp}(\mathcal{K}_{\Exp})$. A similar argument shows that $\sigma_{\Exp}(\mathcal{K}_{\Exp}) \geq \sigma_{\Exp}(\mathcal{K}_{\E_n})$.
\end{proof}

This lemma shows that we can carry out a Hrushovski construction with the class $\mathfrak{C}^{\E_n}$ and predimension $\delta_{\E_n}$. Then we will end up with a strong amalgam $\mathcal{U}_{\E_n}$. The above results imply that $\mathcal{U}_{\Exp}$ and $\mathcal{U}_{\E_n}$ correspond to each other in the $\Exp \leftrightarrow \E_n$ correspondence described at the beginning of this section. Therefore we deduce from Theorem \ref{Adequacy-for-exp} that $\delta_{\E_n}$ is adequate.

\begin{theorem}
Let $\mathcal{F}$ be the countable saturated differentially closed field. Then $\mathcal{U}_{\E_n}$ is isomorphic to $\mathcal{F}_{\E_n}$. In particular, the $\E_n$-reduct of any differentially closed field is elementarily equivalent to $\mathcal{U}_{\E_n}$.
\end{theorem}

We conclude our paper with a final observation. The axiomatisation of $T_{\E_n}$ given in Section \ref{Existential-closedness} is $\forall \exists$. One can also notice that $T_{\E_n}$ is not model complete since otherwise $T_{\Exp}$ would be model complete too, which is not the case (note nevertheless that $T_{\E_n}$ is nearly model complete). In this situation one can apply Theorem 8.1 of \cite{Aslanyan} to conclude that $\D$ is not definable in $T_{\E_n}$. One could also prove this using the fact that $T_{\Exp}$ does not define $\D$ (which can be found in  \cite{Kirby-thesis} and \cite{Aslanyan}).


\bibliographystyle{alpha}      
\bibliography{ref.bib}   

%
%

\end{document}